\documentclass[12pt,reqno]{amsart}
\usepackage{a4wide}
\usepackage[italian,english]{babel}

\usepackage{fancyhdr}
\usepackage{indentfirst}
\usepackage{graphicx}
\usepackage{newlfont} 
\usepackage{dsfont}
\usepackage{amssymb}
\usepackage{amsmath}
\usepackage{latexsym}
\usepackage{amsthm}

\usepackage[colorlinks,allcolors=teal]{hyperref}
\usepackage[capitalize,nameinlink,noabbrev]{cleveref} 
\crefname{page}{page}{pages}

\usepackage{amsthm}

\usepackage[utf8]{inputenc} 

\usepackage{hyperref}
\usepackage{faktor}
\let\origfaktor\faktor
\DeclareRobustCommand{\adaptiveFaktor}[2]{%
    \mathchoice
        {\origfaktor{#1}{#2}}
        {#1/#2}
        {#1/#2}
        {#1/#2}
}
\let\faktor\adaptiveFaktor
\usepackage{mathtools}
\usepackage{thmtools}
\usepackage{tikz}
\usetikzlibrary{matrix,calc}
\usetikzlibrary{cd}
\usepackage{cleveref}

\newtheorem{thm}{Theorem}[section]
\newtheorem{defin}[thm]{Definition}
 
\newtheorem{cor}[thm]{Corollary}
\newtheorem{lemma}[thm]{Lemma}

\theoremstyle{remark}
\newtheorem{rmk}[thm]{Remark}

\newcommand{\N}{\mathbb{N}}
\newcommand{\Q}{\mathbb{Q}}

\newcommand{\Z}{\mathbb{Z}}

\newcommand{\op}{\operatorname}

\newcommand{\mc}[1]{\mathcal{[1]}}

\newcommand{\PP}{\mathbb{P}}

\newcommand{\img}{\operatorname{Im}}

\newcommand{\Mon}{\op{Mon}}

\let\ForAll\forall
\renewcommand{\forall}{\quad \ForAll}

\newcommand{\eps}{\varepsilon}

\newcommand{\dual}[1]{#1\,\check{\vrule height1.3ex width0pt}}

\renewcommand{\div}{\op{div}}
\usepackage{paralist}

\newcommand{\alert}[1]{\par \textcolor{red}{\textbf{\huge{#1}}}\par}
\renewcommand{\alert}[1]{}

\usepackage{csquotes}
\usepackage[backend=biber,url=false,doi=true,isbn=false,eprint=true,
            maxbibnames=99,maxcitenames=99]{biblatex}
\DeclareSourcemap{
  \maps[datatype=bib]{
    \map[overwrite]{
      \step[fieldsource=doi, final]
      \step[fieldset=eprint, null]
      \step[fieldset=eprinttype, null]
      \step[fieldset=eprintclass, null]
    }
  }
}            
\usepackage{todonotes}

\author{Giacomo Nanni}
\address{Dipartimento di Matematica \\
	Universit\`a di Bologna\\
	Piazza di Porta San Donato 5\\
	40127 Bologna, Italy}
\email{giacomo.nanni13@unibo.it}
\address{
    Fakult\"at f\"ur Mathematik\\
    Ruhr-Universit\"at Bochum\\
    IB, Universitätsstraße 150\\
     44801 Bochum, Germany
}
\email{giacomo.nanni@edu.ruhr-uni-bochum.de}

\makeatletter

\@namedef{subjclassname@2020}{\textup{2020} Mathematics Subject Classification}


\title{Lagrangian fibrations on Nikulin-type orbifolds}

\newcommand{\SigmaX}{\Sigma_X}
\newcommand{\SigmaQ}{\ensuremath{\hat{\Sigma}}}

\newcommand{\SigmaY}{{\Sigma_Y}}
\newcommand{\deltaY}{\delta_Y}
\newcommand{\id}{\op{id}}
\newcommand{\Kfour}{$K3^{[2]}$-type }
\newcommand{\LambdaX}{\Lambda_{K3^{[2]}}}
\newcommand{\LambdaFix}{\Lambda_{fix}}
\newcommand{\LambdaY}{\Lambda_{Nik}}
\newcommand{\Mod}{\op{M}}

\newcommand{\Picrk}{\rho_{Pic}}

\begin{document}
\begin{abstract}
    We classify lagrangian fibrations on Nikulin orbifolds, a well studied class of singular irreducible holomorphic symplectic varieties, and prove they verify the SYZ conjecture.
\end{abstract}
\maketitle
\setcounter{tocdepth}{2}
\tableofcontents

\section{Introduction}

Irreducible holomorphic symplectic manifolds are simply connected, compact Kähler manifolds with a unique (up to a scalar factor) holomorphic symplectic form. These manifolds play an important role in the study of complex manifolds with trivial first Chern class, as they constitute one class of irreducible factors in the Beauville--Bogomolov decomposition theorem~\cite{Beauville}. The topology of these manifolds is strongly constrained; in particular their second cohomology group is endowed with a non-degenerate integral quadratic form, the so-called Beauville-Bogomolov-Fujiki form.

A famous conjecture, known as SYZ (see~\cite[Conjecture 1.4]{VerbSYZconj}), predicts that these manifolds admit a lagrangian fibration, meaning a surjective morphism with connected fibers to a positive-dimensional base, if and only if there exists an isotropic nef class in the algebraic part of the second cohomology group.

No general proof of the conjecture is known. However, it has been proved for all known deformation families of IHS manifolds (see~\cite{BayerMacri,MongardiOnorati,MongardiRapagnetta, Yoshioka_BridgelandStability}).

\medskip

Recently, to circumvent the scarcity of examples and to address the natural occurrence of singularities in the minimal model program, singular generalisations of IHS manifolds have been proposed. One possible direction for such generalisation turns back to the original approach by Fujiki~\cite{Fujiki1983}, who considered varieties with at worst finite quotient singularities.

Among the most studied examples are the so-called Nikulin orbifolds (see~\cite{projModels,MenetBBFNik}).
\begin{defin}
    A \textnormal{Nikulin orbifold} is an irreducible symplectic orbifold obtained as the $\Q$-factorial terminalisation of the quotient of a \Kfour fourfold by a symplectic involution.
    Orbifolds deformation equivalent to a Nikulin orbifold are called \textnormal{Nikulin-type orbifolds}.
\end{defin}
More concretely, let $X$ be a $K3^{[2]}$-type fourfold, and let $\sigma\in \op{Aut}(X)$ be a symplectic involution. The quotient of $X$ by $\sigma$
is known to be singular along the disjoint union of a K3 surface $\SigmaQ$ and 28 singular points. The $\Q$-factorial terminalisation $Y$ obtained as the blow-up of $\SigmaQ\subset$
is a Nikulin orbifold.

\begin{rmk}
    When the smooth fourfold $X$ is the Hilbert scheme of two points on a $K3$ surface $S$, and the involution $\sigma$ is induced by a symplectic involution $\sigma_S$ on $S$ as $\sigma=\sigma_S^{[2]}$, the singular locus has a geometric interpretation. It has been shown that a symplectic involution on a $K3$ surface fixes 8 points $p_1,...,p_8\in S$. Then the fixed locus of $\sigma_S^{[2]}$ consists of the K3 surface $\Sigma_X$ coming as the graph of $\sigma_S$ and the $28=\binom{8}{2}$ couples $\{p_i,p_j\}_{i\neq j}$ . More details on the construction can be found in \cite{projModels}.
\end{rmk}

The Beauville-Bogomolov-Fujiki form of Nikulin-type orbifolds has been computed in \cite{MenetBBFNik}, projective models have been provided in \cite{projModels} and the monodromy group has been studied in \cite{MenetMonNik,StevellMonNik,myMonNik}.

The main result of the present work is the following: we provide a classification of lagrangian fibrations on Nikulin-type orbifolds.

\begin{thm}[\Cref{thm:classificationFinal}]
    There are 2 deformation families of lagrangian fibrations on Nikulin-type orbifolds. In particular, let $Y$ be a Nikulin-type orbifold and $\phi: Y\rightarrow B$ a lagrangian fibration, then $\phi$ can be deformed to one of the following examples:
    \begin{itemize}
        \item A fibration induced by a lagrangian fibration defined over a $K3^{[2]}$ fourfold invariant for a symplectic involution.
        \item  A fibration induced by a lagrangian fibration defined over a $K3^{[2]}$ fourfold anti-invariant for a symplectic involution.
    \end{itemize}
\end{thm}

From this classification the SYZ conjecture on Nikulin-type orbifolds follows. This is achieved by the same method used in \cite{OnoratiOrtizSYZsingModuli}, where the conjecture is proven for symplectic varieties deformation equivalent to moduli spaces of stable sheaves.

\subsection{Notation}
Let us fix here the notation that is consistently used throughout the present text.
First, we review the notation we will use for lattices. For any lattice $\Lambda$ we will denote the bilinear pairing as $(\cdot,\cdot)_\Lambda$, dropping the subscript if the context makes it clear. For $v\in \Lambda$, the divisibility of $v$ (meaning the positive generator of the subgroup $(v,\Lambda)$ of $\Z$) will be denoted $\div v$. Also, we will write the square of $v$ as $v^2:=(v,v)$. Finally, for any integer $n$, we will denote by $\Lambda(n)$ the lattice obtained from $\Lambda$ by multiplying the bilinear form by $n$.

Let $X$ be a fourfold of $K3^{[2]}$ type -- that is, deformation equivalent to the Hilbert scheme of two points on a $K3$ surface -- and let $\sigma\in\op{Aut}(X)$ be a symplectic involution. The $K3$ surface in the fixed locus of $\sigma$ is denoted as $\SigmaX$, while the quotient together with its projection is denoted
\[
    \pi\colon  X\longrightarrow \hat X := \faktor{X}{\sigma}.
\]
Finally, let us call $\SigmaQ:=\pi(\SigmaX)\subset \hat X$ the image inside the quotient, and  $Y\xrightarrow{\beta}\hat X$ the blow-up of $\SigmaQ$, with exceptional divisor $\SigmaY$.

Following \cite[Lemma 1.3]{MenetMonNik}, there is a pushforward map \[
    \pi_*:H^*(X,\Z)\rightarrow H^*(\hat X,\Z)\]
such that $\pi_*\pi^*=2\id_{H^*(\hat X,\Z)}$.

The cohomology groups $H^2(X,\Z)$ and $H^2(Y,\Z)$ will always be considered equipped with the Beauville-Bogomolov-Fujiki quadratic forms. These are respectively isomorphic to the abstract lattices:
\begin{align*}
    \LambdaX & :=  U^3\oplus E_8(-1)^2\oplus \langle-2\rangle \\
    \LambdaY & := U^3(2)\oplus E_8(-1)\oplus \langle
    -2\rangle^2.
\end{align*}
Up to choosing the identification $ H^2(X,\Z)\cong\LambdaX$ one can assume $\sigma^*$ acts by exchanging the two $E_8(-1)$ factors.
Therefore, the $\sigma^*$-invariant part $H^2(X,\Z)^\sigma\subset H^2(X,\Z)$  is isomorphic to $$\LambdaFix:=U^3\oplus E_8(-2)\oplus\langle-2\rangle.$$ The map \begin{equation}\label{def:latticeEmbeddingFixedPart}
    \eta:=\beta^*\circ\pi_*:H^2(X,\Z)^\sigma(2)\rightarrow H^2(Y,\Z)
\end{equation} is a non-primitive isometric embedding of lattices, which at the level of abstract lattices can be written as $u\oplus e\oplus\alpha\mapsto u\oplus 2e\oplus \alpha\oplus\alpha$ for $u\in U^3,e\in E_8$ and $\alpha\in \langle-2\rangle$.

\medskip

\noindent \textbf{Acknowledgments.}The author thanks his supervisors Christian Lehn and Giovanni Mongardi, for introducing him to the problem and for the useful discussions. A special thanks to Alessandro Frassineti, for his help in understanding the computations of \cref{lemma:generalPicMT}, and to Claudio Onorati for the discussions around \cref{thm:classificationFinal}.
The author is part of INdAM research group “GNSAGA” and was partially supported by it. The author acknowledges the support of the European Union - NextGenerationEU under the National Recovery and Resilience Plan (PNRR) - Mission 4 Education and research - Component 2 From research to business - Investment 1.1 Notice Prin 2022 - DD N. 104 del 2/2/2022, from title "Symplectic varieties: their interplay with Fano manifolds and derived categories", proposal code 2022PEKYBJ – CUP J53D23003840006. The author was also supported by the DFG Project-ID 550535392 and DFG Project-ID 530132094.

\section{Preliminaries}
\subsection{Examples of fibrations}\label{sec: examples}
We give now some examples of Lagrangian fibrations on orbifolds of Nikulin type. These will be the prototypical examples for the classification in \cref{thm:classificationFinal}.

\subsubsection{On the Fermat quartic}
Consider the Fermat quartic $S=V(\sum_{i=0}^3x_i^4)\subset \PP^3$.
There is a projection from $S$ over a smooth quadric surface $Q$ which can be written explicitly as $[x_0,...,x_3]\mapsto[x_0^2,...,x_3^2]$. As $Q\cong \PP^1\times \PP^1$, by composing with one of these two projections we get an elliptic fibration $\eps:S\rightarrow \PP^1$. Explicitly, for some coordinate choice, it can be written as \[
    \eps([x_0,x_2,x_3,x_4])=[x_0^2+\sqrt{-1}x_1^2,x_2^2+\sqrt{-1}x_3^2].\]

Let us consider the involutions \[
    \sigma_+([x_0,x_1,x_2,x_3]):=[-x_0,-x_1,x_2,x_3] \quad \text{and} \quad \sigma_-(([x_0,x_1,x_2,x_3]))=[x_2,x_3,x_1,x_0].\]
By a direct computation of the pullback of the symplectic form, one checks that both these involutions are symplectic.
Indeed, $\eps\circ \sigma_+=\eps$, while $\eps\circ\sigma_- =\sigma_-' \circ\eps$, where $\sigma_-'\in PGL(2)$ is the non trivial involution of $\PP^1$, $\sigma_-'([y_0,y_1]):=[y_1,y_0]$.

We can induce involutions $\sigma_\pm^{[2]}\in \op{Aut}_s(S^{[2]})$ and a lagrangian fibration \[
    \eps^{[2]}:S^{[2]}\rightarrow\op{Sym}^2(\PP^1)\cong\PP^2.\]
The K3 surface $\Sigma_\pm\subset S^{[2]}$ fixed by $\sigma_\pm^{[2]}$ is the graph of $\sigma_\pm$, namely
\[
    \Sigma_\pm=\{\{p,\sigma_\pm(p)\}\in S^{[2]}, p\in S\}.
\]
The image of $\Sigma_+$ through $\eps^{[2]}$ corresponds to double points $\{\eps(p),\eps(p)\}\in\op{Sym}^2\PP^1$. Under the identification
\[
    \PP^2\stackrel{\cong}{\longrightarrow}\op{Sym}^2\PP^1, \quad [a_0,a_1,a_2] \longmapsto V(\sum a_iy_0^iy_1^{2-i}),\] such double points correspond to the zero locus of $a_0a_2-2a_1^2$ inside of $\PP^2$. In particular, this shows that $\Sigma_+$ does not intersect the generic fiber. An analogous identification is carried out for $\Sigma_-$, for which $\eps^{[2]}(\Sigma_-)$ corresponds to the the line $V(t_0-t_2)$.

Denote $\sigma_+'=\id_{\PP^2}$; then the fibration $\eps^{[2]}$ induces maps
\[
    \faktor{S^{[2]}}{\sigma_\pm}\rightarrow \faktor{\PP^2}{\sigma_\pm '}.\] For each one of them, blowing up the K3 surface $\Sigma_\pm$ does not affect the generic fiber. Thus, in both cases, we get a lagrangian fibration on the Nikulin orbifold \[
    \phi_\pm:Y_\pm=\op{Bl}_{\Sigma_\pm}\faktor{S^{[2]}}{\sigma_\pm}\rightarrow \faktor{\PP^2}{\sigma_\pm '}\] respectively one with source $Y_+$ and base $\PP^2$ and the other with source $Y_-$ and base $\PP^2/\tau$.

\begin{rmk}\label{rmk:polarisationTypesEquivariant}
    As the generic fiber of $\phi_-$ is isomorphic to that of $\eps^{[2]}$, it is principally polarised by \cite[Theorem 6.1]{Wieneck_polarizationTypes}, while the generic fiber of $\phi_+$ is a quotient of a principally polarised abelian surface by an involution without fixed points, hence it is polarised of type $(1,2)$ by \cite[Section 1]{Barth_AbelianSurfaces12Polarization}.
\end{rmk}

\subsubsection{The Markushevich-Tikhomirov Prym}\label{ex:MTConstruction}
In \cite{MarkTikPrym}, a Nikulin-type orbifold is constructed as a Prym variety. Their example comes with a natural lagrangian fibration of polarisation type (1,2) (see \cite[Theorem 3.4]{MarkTikPrym}). We recall the construction.
\newcommand{\vv}{\mathbf{v}}
Let
\[
    \mu:P\xrightarrow[\mathcal{O}(4)]{2:1}\PP^2\]
be a double cover ramified over a quartic $Q_0\in |\mathcal{O}(4)|$. Let
\[
    \rho:S\xrightarrow[-2K_P]{2:1}P
\]
be a double cover ramified over a section $\Delta_0\in |-2K_P|$, namely linearly equivalent to twice the canonical bundle $K_P$ of $P$. We denote $Q,\Delta$ respectively the preimage of $Q_0,\Delta_0$ via $\mu,\rho$. We denote $\tau$ the cover involution on $S$ and $\alpha:=\mu\circ \rho$.

For each line $\ell\subset \PP^2$ we have a curve in $P$ defined as $D_\ell:=\mu^{-1}(\ell)$ and one in $S$ as $C_\ell:=\rho^{-1}(D_\ell)$. The double cover $C_\ell\rightarrow D_\ell$ has cover involution $\tau_\ell:=\tau_{|C_\ell}$ and defines a Prym variety $\op{Prym}(C_\ell,\tau_\ell)$ as the fixed locus in the Jacobian $J(C_\ell)$ of the involution $-\tau_\ell^*$. Markushevich and Tikhomirov produce a compactification $\cal P$ of this family, realised inside the moduli space $M:=\Mod_H(0,H,-2)$ of semistable sheaves on $S$ with respect to the ample class $H:=\alpha^*\mathcal{O}_{\PP^2}(1)$ and Mukai vector $\vv=(0,H,-2)$ (see \cite[Definition 3.3]{MarkTikPrym}). The natural support map ${\cal P}\rightarrow|H|$ is a lagrangian fibration, and moreover $\cal P$ is birational to a Nikulin orbifold \cite[Theorem 3.4, Corollary 5.7]{MarkTikPrym}. We recall here the description of the birational map.

Consider the rational map $\tilde\Phi:S^{[2]}\dashrightarrow {\cal P}$ sending $s,t\mapsto [s+\tau(s)]+[t+\tau(t)]$ where we assume $\pi(s)\neq\pi(t)$ and the brackets denote the associated divisor class in $Prym(C_\ell,\tau_\ell)$ with $\ell=<\alpha(s),\alpha(t)>\subset \PP^2$ the line between $\pi(s)$ and $\pi(t)$. This map is invariant under composition with the Beauville involution $\iota_0\in \op{Aut}(S^{[2]})$ and with the involution $\tau^{[2]}$ induced on $S^{[2]}$ by $\tau$. Quotienting by $\sigma=\beta\circ \tau^{[2]}$ gives a birational map
\[
    \Phi \colon \faktor{S^{[2]}}{\sigma}\dashrightarrow {\cal P},\]
and in turn the birational map $\Phi\circ \beta:Y\dashrightarrow \cal P$ obtained by composing with the blow up of the K3 surface fixed by $\sigma$. Since both $Y$ and $\cal P$ are $\Q$-factorial and terminal they are deformation equivalent \cite[Theorem 4.9]{BakkerLehn} so $\cal P$ is deformation equivalent to a Nikulin orbifold.

\subsection{Monodromy}
Let $Y$ be an orbifold of Nikulin type (or more generally an irreducible holomorphic symplectic varieties, see \cite[Definition 4.1]{BakkerLehn}).
\begin{defin}
    A locally trivial deformation of $Y$ is a proper flat morphism  $\cal Y\rightarrow T$  such that for some point $0\in T$ its fiber $\cal Y_0\cong Y$ and  for all $t\in T, y\in \cal Y_p$ there exists an open sets $\cal U\subset \cal Y, V\subset T$ such that $\cal U\cong \cal Y_p\times U$.
\end{defin}
\begin{defin}
    Let  $\cal Y\rightarrow T$ be a locally trivial deformation of $Y$ and $t\in T$.
    A parallel transport operator $H^2(Y,Z)\rightarrow H^2(Y_t,\Z)$ is a map $P_\gamma:H^2(Y,\Z)\rightarrow H^2(Y_t,\Z)$ induced by parallel transport along a path $\gamma$ in $T$ connecting $0$ to $t$.
\end{defin}
If the base $T$ is not simply connected, the fundamental group $\pi_1(T,0)$ may act non-trivially on the cohomology of $Y$: any path $\gamma\in \pi_1(T,0)$ induces in particular a map $H^2(Y,\Z)\rightarrow H^2(Y,\Z)$.
Parallel transport operators are isometries for the Beauville-Bogomolov-Fujiki quadratic form, forming a subgroup of the orthogonal group $O(H^2(Y,\Z))$.
\begin{defin}
    The monodromy group $\Mon^2(Y)$ is the subgroup of isometries of $H^2(Y,\Z)$ formed by parallel transport operators.
\end{defin}
Given a marking, which means an isometry $\eta:H^2(Y,\Z)\rightarrow \LambdaY$, this induces an action of $\Mon^2(Y)$ on $\LambdaY$ by conjugation via $\eta$. We denote $\Mon^2_\eta (\LambdaY)$ the image of this action in the group of isometries of $\LambdaY$. This action depends (up to conjugation by an isometry of $\LambdaY$) on the chosen marking. However, recent works (as \cite{StevellMonNik,myMonNik}) have given a numerical characterisation of $\Mon^2 _\eta \LambdaY$ showing that it does not depend on $\eta$. We can therefore drop the subscript from the notation and write $\Mon^2 (\LambdaY)$.

\subsection{Isotropic classes}\label{sec:isotropicClasses}
Let $Y$ be a Nikulin-type orbifold and $\phi:Y\rightarrow B$ a lagrangian fibration. By \cite[Theorem 3]{SchwaldRelationsAndFibrationsIHSV}, the Picard group $\op{Pic}(B)$ is of rank 1, generated by an ample class $H\in \op{Pic}(B)$. The pullback $l:=\phi^*c_1(H)\in H^2(Y,\Z)$ is called the \emph{class} of the lagrangian fibration. From the Fujiki relations, it follows that $l^2=0$ where the product is taken in the sense of the Beauville-Bogomolov-Fujiki form of $Y$.

In \cite[Theorem 6.15]{MenetMonNik} some monodromy operators were produced and representatives were found for each monodromy orbit.
For the convenience of the reader, we report here their description.
Let us denote $L_i\in U(2)$ a primitive class of square $L_i^2=4i, i\in \N$, $e_i\in E_8(-1)$ a primitive class of square $e_i^2=-2i$ for $i=1,2$ and $\gamma_1,\gamma_2$ two primitive generators for the two $\langle -2 \rangle$  factors. Denote $\deltaY:={\gamma_1+\gamma_2}$ and $\SigmaY:={\gamma_1-\gamma_2}$.

\begin{thm}{\cite[Theorem 6.15]{MenetMonNik}} \label{thm:9monorb-M'}
    Let $v\in \LambdaY$ be a primitive non-zero element. Denote by $v_{E_8}$  the projection of $v$ to the $E_8(-1)$-part of the lattice, and let ${\bar{v}}_{E_8}$ be its image in the $\Z/4\Z$-module $E_8(-1)/4E_8(-1)$.
    Then there exists a monodromy operator
    $f\in \Mon^2(\LambdaY)$ such that
    $$ f(v)=\left\{
        \begin{array}{l}
            \textrm{If $v$ satisfies $(*)$:}                                                        \\
            \hspace{0.5 em}
            1)  \hspace{1em}L_{i}  \hspace{7.5em}\textrm{with\ }\div(v)=2 \textrm{\ and\ } q(v)=4i. \\
            \textrm{Otherwise,
            }                                                                                       \\
            \begin{array}{lllll}
                2) & 2L_{i} - \deltaY                           & \, \textrm{if\ } \div(v)=2,  & q(v)=16i-4,  & \textrm{and\ }  {\bar{v}}_{E_8}=0           \\
                3) & 2L_{i+1} + 2e_2 -\deltaY                   & \, \textrm{if\ } \div(v)=2,  & q(v)=16i-4,  & \textrm{and\ } {\bar{v}}_{E_8}\neq 0        \\
                4) & L_i - \gamma_1                             & \, \textrm{if\ } \div (v)=2, & q(v)=4i-2,   & \textrm{and\ } {\bar{v}}_{E_8}=0            \\
                5) & L_{i+1} + e_2-\gamma_1                     & \, \textrm{if\ } \div (v)=1, & q(v)=4i-2,   & \textrm{and\ }
                q(v_{E_8})\equiv 0 \pmod{4}                                                                                                                 \\
                6) & L_{i} + e_1                                & \, \textrm{if\ } \div (v)=1, & q(v)=4i-2,   & \textrm{and\ } q(v_{E_8})\equiv 2 \pmod{4}  \\
                7) & 2L_{i} + 2e_1 - \deltaY                    & \, \textrm{if\ } \div(v)=2,  & q(v)=16i-12, & \textrm{and\ } {\bar{v}}_{E_8}\neq 0        \\
                8) & L_i + e_1- \gamma_1                        & \, \textrm{if\ } \div(v)=1,  & q(v)=4i-4,
                   & \textrm{and\ } q(v_{E_8})\equiv 2 \pmod{4}                                                                                             \\
                9) & L_{i+1} + e_2                              & \, \textrm{if\ } \div(v)=1,  & q(v)=4i,     & \textrm{and\ } q(v_{E_8})\equiv 0 \pmod{4}. \\
            \end{array}
        \end{array}
        \right. $$
    The vector $v$ satisfy $(*)$ if
    \begin{compactenum}
        \item The restriction $v_{U^3(2)}$ of $v$ to $U^3(2)$ is not divisible by 2,
        \item the restriction $v_{E_8}$ of $v$ to $E_8(-1)$ is divisible by 2, and
        \item the restriction $v_{(-2)\oplus (-2)}$ to  $\langle \frac{\deltaY+ \SigmaY}{2} ,  \frac{\deltaY-
                \SigmaY}{2}\rangle$ is contained in the sublattice $\langle \deltaY, \SigmaY\rangle$.
    \end{compactenum}

\end{thm}
\begin{rmk}
    Notice that the result only asserts that a vector $v$ satisfying the numerical condition in one of the bullet point is monodromy equivalent to the given class. The vice versa is not true as we see in the next corollary.
\end{rmk}
\begin{cor}\label{cor:numberOrbits}
    There are two monodromy orbits for primitive isotropic classes. In the notation of \cref{thm:9monorb-M'}, they are generated by $L_0,L_1+e_2$.
    \begin{proof}
        From \cref{thm:9monorb-M'} there are at most three monodromy orbits for primitive isotropic vectors: cases 1 and 9 with $i=0$ and case 8 with $i=1$. The case 1 has different divisibility than the others, so there is no isometry (and in particular no monodromy operator) sending $L_0$ in $L_1+e_2$ or $L_1+e_1+\gamma_1$. However, cases 8 and 9 are monodromy equivalent. To prove it, assume that $\eps_1,...,\eps_8$ is the basis of simple roots of $E_8(-1)$ so that the pairing has associated matrix the Cartan matrix. Then we can assume $e_2=\eps_1+\eps_3$ and $e_w=\eps_4+\eps_6$, so that $e_2^2=e_w^2=-4$ and $(e_2,e_w)=1$. Denote $w:=L_1+e_w+\gamma_1$: the vector $w$ is primitive, of divisibility 1 and square $L_1^2+e_w^2+\gamma_1^2=4-4-2=-2$.  The reflection $R_v:\LambdaY\rightarrow\LambdaY$ defined as $x\mapsto x+(x,v)v$ is a reflection around a vector of negative square and therefore a monodromy operator (see also \cite[Corollary 3.1 ]{StevellMonNik}). We can check that $(L_1+e_2,w)=L_1^2+(e_2,e_w)=5$ and so $R_v(L_1+e_2)=L_1+e_2+5w$. The component on the $E_8(-1)$ summand is $R_v(L_1+e_2)_{E_8}=e_2+5e_w$ which has square $(e_2+5e_w)^2=e_2^2+25e_w^2+2(e_2,e_1)\cong 2 \mod 4$. Therefore, there exists a monodromy operator $f\in \Mon(\LambdaY)$ such that $f\circ R_w(L_1+e_2)=L_1+e_1+\gamma_1$.\endproof
    \end{proof}
\end{cor}
As a consequence of this description, we get an upper bound on the monodromy orbits of classes of lagrangian fibrations.
\begin{lemma}\label{lemma:max3Orbits}
    In each integral ray $\Z_{>0}\hat l$, with $\hat l$ a primitive isotropic class, there is at most one class associated to a lagrangian fibration. In particular, there are at most 2 monodromy orbits of lagrangian fibrations on Nikulin-type orbifolds.
    \proof
    Suppose there exist two lagrangian fibrations $\phi_i:Y\rightarrow B_i$ of a Nikulin-type orbifold $Y$ with associated classes respectively $l_i=k_i\hat l$. As $\mathcal{O}_{B_1}(1)$ is ample, some multiple $\mathcal{O}_{B_1}(t_1)$ defines an embedding $\psi_1:B_1\hookrightarrow \PP^{n_1}$. The composition with the Veronese map ${\mathcal{V}}_{t_2k_2}:\PP^{n_1}\hookrightarrow \PP^{N_1}$ of degree $t_2k_2$ satisfies
    \[
        \phi_1^*\psi_1^*\mathcal{V}_{t2k2}^*\mathcal O_{\PP^{N_1}}(1)=k_1k_2t_1t_2 \hat l.\]
    The same holding exchanging the roles of 1 and 2 shows, by \cite[II.7.1]{Hartshorne}, that the two maps  $\mathcal{V}_{t2k2}\circ \psi_1\circ \phi_1,\mathcal{V}_{t1k1}\circ \psi_2\circ \phi_2$ only differ by an automorphism of the target. Therefore, as the Veronese and the $\psi_i$ are embeddings, we get an automorphism $\gamma:B_1\cong B_2$ such that $\gamma\circ \phi_1=\phi_2$. This implies that $k_1=k_2$.

    From \cref{cor:numberOrbits}, that there are only 2 orbits of primitive isotropic vectors.
    \endproof
\end{lemma}
\begin{rmk}
    The first part of the proof of \cref{lemma:max3Orbits} does not rely on anything specific about Nikulin-type orbifolds, and it works the same for more general irreducible holomorphic symplectic varieties.
\end{rmk}

\begin{rmk}\label{rmk:differentNumericalProperties}
    Notice that the two classes $L_0,L_1+e_2$ have different numerical properties: $\div L_0=1$ while $\div L_1+e_2=2$. So there is no isometry of $\LambdaY$ sending a vector in the first monodromy orbit in the second one.
\end{rmk}
We give the following definition:
\begin{defin}\label{def:types}
    Let $Y$ be an orbifold of Nikulin type, $\psi:Y\rightarrow B$ a lagrangian fibration of class $l_\psi=\psi^*c_1(\mathcal{O}_B(1))=k\hat l$ with $k\in \Z_{>0}$ and $\hat l\in H^2(Y,\Z)$ primitive. Then $\psi$ is said to be of type A (resp. B) if there exists a marking $\eta:H^2(Y,\Z)\rightarrow \LambdaY$ such that $\eta(\hat l)=L_1+e_2$ (resp. $L_0$).
\end{defin}

Notice that by \cref{lemma:max3Orbits} any lagrangian fibration satisfies one of these two conditions (type A,B ) and by \cref{rmk:differentNumericalProperties} it satisfies exactly one.
\begin{rmk}\label{rmk:anyMarking}
    Notice that by \cref{rmk:differentNumericalProperties}, as two markings differ by an isometry of $\LambdaY$, any marking on a Nikulin-type orbifold with a lagrangian fibration of type A (or B) can be assumed to satisfy the requirement of \cref{def:types} up to composing with some monodromy operator $P_{\omega_1}\in Mon^2(Y)$. Indeed, given any marking $\eta$ on $Y$, where $\psi:Y\rightarrow B$ is a fibration of type A (similarly the same holds for type B) from \cref{cor:numberOrbits} there exists a monodromy operator $g\in \Mon^2(\LambdaY)$ such that $g\eta(l_\psi)\in \{L_0,L_1+e_2\}$. Let $\hat\eta$ the marking on $Y$ such that $\hat \eta (l_\psi)=L_0$. Then $g\eta\hat\eta^{-1}$ is an isometry sending $L_0$ to $g\eta(l_\psi)$ which therefore must be $L_0$. Using the identification induced by $\eta$ between $\Mon(\LambdaY)$ and $\Mon(Y)$ we get that there exists a monodromy operator $P_{\omega_1}\in \Mon^2(Y)$ such that $\eta P_\omega=g\eta$ is a marking on $Y$ satisfying the condition in the definition of type A: $\eta P_\omega(l_\psi)=L_0$.
\end{rmk}

\section{Classification of fibrations}

In this section we obtain the main result, \cref{thm:classificationFinal}, which provides a classification of Lagrangian fibrations on Nikulin orbifolds up to monodromy operators. For each one of the examples exposed in \cref{sec: examples} we compute the corresponding monodromy orbit containing the class of the lagrangian fibration. First in \cref{sec:equivariantFibrations} we discuss the two cases arising from lagrangian fibrations on the smooth $K3^{[2]}$ fourfold. Then in \cref{ssec:MTsystem} we show that the Markushevich-Tikhomirov system is of type A.
\subsection{Equivariant fibrations on the smooth fourfold}\label{sec:equivariantFibrations}

Let $X$ be a smooth \Kfour fourfold, $\sigma\in \op{Aut}(X)$ a symplectic involution and $\psi:X\rightarrow \PP^2$ a lagrangian fibration with class $l_X:=\psi^*\mathcal{O}_{\PP^2}(1)$ such that $\sigma^*l_X=l_X$. By the $\sigma$-invariance assumption on $l_X$, there is an induced involution (which we denote $\sigma'$) on $\PP^2$ such that $\phi\circ\sigma=\sigma'\circ\phi$. By passing to quotients, $\psi$ induces the map $\hat\psi:\faktor{X}{\sigma}\rightarrow \faktor{\PP^2}{\sigma'}$ and therefore a lagrangian fibration $\phi:=\beta\circ\hat\psi:Y\rightarrow \faktor{\PP^2}{\sigma'}$ on the Nikulin orbifold $Y$ constructed as $\Q$-factorial terminalisation of $\faktor{X}{\sigma}$. The commutative diagram below sums up the construction.
\begin{center}
    \begin{tikzcd}
        \phantom{a}&X\arrow{d}{\pi}\arrow{r}{\psi}&\PP^2\arrow{d}{\hat\pi}\\
        Y\arrow{r}{\beta}&\faktor{X}{\sigma}\arrow{r}{\hat\psi}&\faktor{\PP^2}{\sigma'}
    \end{tikzcd}
\end{center}

Since $\sigma$ is an involution over a projective plane it can either be $\sigma'=\id_{\PP^2}$, or it can be
\[
    \sigma'([a_0,a_1,a_2])=[-a_0,a_1,a_2], [a_0,a_1,a_2]\in\PP^2.\] In the following, we describe the monodromy orbit for the class $l_Y$ in the two cases.
We say that the fibration $\psi$ is invariant if $\sigma'=\id_{\PP^2}$ and anti-invariant otherwise.

\subsubsection{Trivial action on the base: type A}\label{ssec:invariant}
We first address the case of invariant fibrations -- that is lagrangian fibrations $\psi:X\rightarrow \PP^2$ such that $\psi\circ\sigma=\psi$.
\begin{lemma}\label{lemma:invariant}
    In the notation from before, assume $\sigma'=\id_{\PP^2}$. Then $2l_Y=\eta(l_X)$ and the fibration is of type A.
    \proof
    By definition, we get $l_\phi=\beta^*\hat\psi^*\hat H$ where $\hat H$ is the generator of the Picard group of $\faktor{\PP^2}{\sigma'}$. By \cite[Proposition 1.3]{MenetBBFNik}, we have $\pi_*\pi^*=2$ so  we can write $$2l_\phi=\beta^*\pi_*\pi^*\hat\psi^*\hat H=\eta(\psi^*\mathcal{O}(k))=k\eta(l_\psi).$$
    As $\faktor{\PP^2}{\sigma'}=\PP^2$ and $\hat \pi=\id_{\PP^2}$,  we have
    $$2l_Y=\beta^*2\hat\psi^*\mathcal{O}(1)=\beta^*\pi_*\pi^*\hat\psi^*\mathcal{O}(1)=\eta(l_\psi).$$
    Since $l_\psi\in H^2(X,\Z)$ is primitive, its image $\eta(l_\psi)=2l_\phi$ is primitive inside the sublattice $\img \eta\subset H^2(Y,\Z)$. Since the $\img \eta$ has index 2 in its saturation, $l_\phi$ is primitive. However, notice that $l_\phi\not\in\img \eta$ as $2l_\phi$ is primitive in $\img \eta$. Therefore, the $E_8(-1)$ component must be an odd multiple of a primitive vector, implying that $\div l_\phi=1$. This concludes by \cref{cor:numberOrbits}.
    \endproof
\end{lemma}

\subsubsection{Non-trivial action on the base: type B} \label{ssec:antiInvariant}
We consider now the case of the fibrations induced by anti-invariant fibrations of the smooth $K3^{[2]}$ fourfold, meaning fibrations $\psi:X\rightarrow \PP^2$ such that $\psi\circ\sigma=\sigma'\circ\psi$ where $\sigma'$ is a non-trivial linear involution on $\PP^2$.
\begin{lemma}\label{lemma:antiInv}
    In the notation from before, assume $\sigma'\neq\id_{\PP^2}$. Then $l_Y=\eta(l)$ and the fibration is of type B.
    \proof
    The quotient $\faktor{\PP^2}{\sigma'}$ is isomorphic to $\PP(1,1,2)$. Let $\mathcal{O}_{112}(2)$ be the ample generator of the Picard group of $\PP(1,1,2)$ (see \cite[Section 1.4]{DolgachevWeighted} for more details on weighted projective spaces).
    Let $v:\PP^2\rightarrow \PP^5$ be the Veronese map and  $\hat v:\PP(1,1,2)\rightarrow \PP^3$ the map induced by $\mathcal{O}_{112}(2)$. We have the following commutative diagram.
    \begin{center}
        \begin{tikzcd}
            &X\arrow{d}{\pi}\arrow{r}{\psi}&\PP^2\arrow{r}{v}\arrow{d}{}&V\arrow[hookrightarrow]{r}{}{\subset}\arrow{d}{p}&\PP^5\arrow[dashed]{d}{}\\
            Y\arrow{r}{\beta}&\faktor{X}{\sigma}\arrow{r}{\hat\psi}&\PP(1,1,2)\arrow{r}{\hat v}&\hat V\arrow[hookrightarrow]{r}{\subset}&\PP^3
        \end{tikzcd}
    \end{center}
    where $p:V\rightarrow\hat V$ is the restriction of the projection from a line $\PP^5\dashrightarrow\PP^3$.

    Observe that $\mathcal{O}_{112}(2)=\hat v^*\mathcal{O}_{\PP^3}(1)_{|\hat V}$. Therefore, similarly to \cref{lemma:invariant}, we write: $$2l_Y=2\beta^*\hat\psi^*\hat v^* \mathcal{O}_{\PP^3}(1)=\beta^*\pi_*\pi^*\hat\psi^*\hat v^* \mathcal{O}_{\PP^3}(1)_{|\hat V}=\eta \psi^*v^*p^*\mathcal{O}_{\PP^3}(1)_{|\hat V}.$$ But since $p^*\mathcal{O}_{\PP^3}(1)_{|\hat V}=\mathcal{O}_{\PP^5}(1)_{|V}$ we get $2l_Y=\eta \psi^*\mathcal{O}_{\PP^2}(2)=\eta (2 l)$ implying $l_Y=\eta (l)$. Since $l_Y\in \img \eta$ it has even divisibility. \endproof
\end{lemma}

\subsection{The Markushevich-Tikhomirov system}\label{ssec:MTsystem}

We now classify the type of the Markushevich-Tikhomirov fibration. We signal for the reader interested in \cref{thm:classificationFinal,cor:SYZ}, that this section is not necessary for the main results. Its main result (\cref{lemma:classMT}) could even be deduced directly from  \cref{cor:orbitsPolarisations}, as it shows that the type of the fibration (in the meaning of \cref{def:types}) is determined a posteriori by the polarisation type of the general fibre. However, we think that this direct approach is interesting in its own.

We start by a numerical consideration. We keep the notation from \cref{thm:9monorb-M'}.

\begin{lemma}\label{lemma:discriminateThirdOrbit}
    Let $v\in \LambdaY$ a primitive isotropic vector. If $(v,\Sigma_Y)\equiv 2 \op{mod} 4$ then $v$ is in the monodromy orbit of $L_1+e_2$

    \proof
    Suppose $v$ is in the monodromy orbit of $L_0$. Since $\div v=2$, we can write  $v=u+2e+k\gamma_1+m\gamma_2$ for $u\in U^3(2),e\in E_8(-1),k,m\in\Z$.
    Since $v^2=0$ we have $$0=u^2+4e^2+-2k^2-2m^2\equiv 2(k^ 2+m^ 2) \mod 4.$$
    On the other hand, $$(v,\Sigma_Y)=k(\gamma_1,\Sigma_Y)+m(\gamma_2,\Sigma_Y)=-2k+2m=-2(k-m)\equiv_4 -2k+2m\equiv_4 -2(k+m)\equiv_4 0.$$
    \endproof
\end{lemma}
Before entering the main result (\cref{lemma:classMT}) we recall some known facts from \cite{MarkTikPrym} and we check some preliminary identities.

Recall from \cref{ex:MTConstruction} that there is a birational map $Y\dashrightarrow \cal P$ between the $\Q$-factorial terminalisation of $\faktor{S^{[2]}}{\sigma}$ and the Prym variety associated to the double cover $S\rightarrow P$. This birational map is a flop of a projective plane, as stated in \cite[Corollary 5.7]{MarkTikPrym}. This flop is induced by a rational $\sigma$-invariant map $S^{[2]}\dashrightarrow \cal P$. A resolution of indeterminacy for this last map is obtained via successive blow ups. It can be described as the incidence variety $$N:=\{(\xi,\ell)\in S^{[2]}\times \dual{\PP^2}, \alpha_*\xi\subset \ell\}.$$
The following diagram fixes the notation for this section
\begin{center}
    \begin{tikzcd}
        N\arrow[r,"\tilde \theta_Y"]\arrow[d,"\hat\pi"]\arrow[rrr,bend left=30,"\tilde\theta_\cal P"]&M\arrow[r,"\tilde\beta"]\arrow[d,"\tilde \pi"]&X=S^{[2]}\arrow[r,dashed,"\tilde\Phi"]\arrow[d,"\pi"]&\cal P\arrow[r,"\phi"]&\dual{\PP^2}\\
        \faktor{N}{\sigma_N}\arrow[r,"\theta_Y"]\arrow["\theta_\cal P",rrru,bend right=50]&\faktor{M}{\sigma_M}=Y\arrow[r,"\beta"]&\faktor{X}{\sigma}\arrow[ru, dashed,"\Phi"]
    \end{tikzcd}
\end{center}
We will denote $l_\cal P:=\phi^*\mathcal{O}_{\dual{\PP^2}}(1)$ the class of the fibration $\phi$.

For a rational map, one can define the pullback of divisors through a resolution of indeterminacy. In particular, we define $\tilde\Phi^*:=\beta_*\tilde{\theta_Y}_*\tilde\theta_\cal P^*$ and similarly $\Phi^*:=\beta_*{\theta_Y}_*{\theta_\cal P}^*$ and $(\Phi\circ \beta)^*:={\theta_Y}_*{\theta_\cal P}^*$. Notice that, in general, this pullback definition is not functorial.

Let $l_X:=\tilde \Phi^*l_\cal P$. If pullback were functorial, it would be clear that it is $\sigma$-invariant by commutativity of the diagram. Fortunately, this is still true.
\begin{lemma}\label{lemma:lXinvariance}
    The class $l_X$ is $\sigma$-invariant.
    \proof
    Since $\sigma^*={\sigma^{-1}}_*=\sigma_*$,  the explicit computation
    \[
        \sigma^* l_X
        = (\sigma_* \tilde{\beta}_* \tilde{\theta}_{Y*} \tilde{\theta}_{\mathcal{P}}^*) l_Y
        = (\tilde{\beta}_* \tilde{\theta}_{Y*} \sigma_{N*} \tilde{\theta}_{\mathcal{P}}^*) l_Y
        = (\tilde{\beta}_* \tilde{\theta}_{Y*} \sigma_N^* \tilde{\theta}_{\mathcal{P}}^*) l_Y
        = (\tilde{\beta}_* \tilde{\theta}_{Y*} (\tilde{\theta}_{\mathcal{P}} \circ \sigma_N)^*) l_Y
        = l_X
    \]
    proves the statement.
    \endproof
\end{lemma}
Let $l_Y:=(\Phi\circ\beta)^*l_\cal P$ the pullback of the class of the lagrangian fibration $\phi:\cal P\rightarrow\dual {\PP^2}$. The two classes $l_Y,l_X$ are related by the following:
\begin{lemma}\label{lemma:lXexplicit}
    In the previous notation, $$\pi^*\beta_*l_Y=l_X.$$
    \proof
    First, one can observe that
    $$\pi^*\beta_*l_Y=\pi^*\beta_*{\theta_Y}_*\theta_{\cal P}^*l_\cal P=\pi^*\Phi^*l_\cal P.$$
    Again, we cannot conclude directly by functoriality. We can however compute $$\Phi^*=\beta_*{\theta_Y}_*\theta_{\cal P}^*=\frac{1}{2}\beta_*{\theta_Y}_*\hat \pi_*\hat\pi^*\theta_{\cal P}^*=\frac{1}{2}(\pi\circ\tilde\beta\circ\tilde\theta_Y)_*\tilde\theta_{\cal P}^*=\frac{1}{2}\pi_*{\tilde\Phi}^*.$$
    Which gives
    $$\pi^*\Phi^*l_{\cal P}=\frac{1}{2}\pi^*\pi_*l_X.$$
    By \cite[Corollary 1.4]{MenetBBFNik}, $\pi^*\pi_*l_X=l_X+\sigma^*l_X$. Since $l_X$ is $\sigma$-invariant by \cref{lemma:lXinvariance}, one has $$\pi^*\beta_*l_Y=\frac{1}{2}\pi^*\pi_* l_X=l_X$$ concluding. \endproof
\end{lemma}
\begin{cor}\label{cor:lYExplicit}
    The class $l_Y$ can be written as
    $$2l_Y=\beta^*\pi_*l_X+k\Sigma_Y$$
    for some $k\in \Z$.
    \proof
    As $\beta^*\pi_*(H^ 2(S^{[2]},\Z)^\sigma)\oplus \Z \Sigma_Y$ is an index 2 sublattice, there exists a $\sigma$-invariant class $L\in H^ 2(S^{[2]},\Z)^\sigma$ and $k\in \Z$ such that:
    $$2l_Y=\beta^*\pi_*L+k\Sigma_Y.$$
    By \cref{lemma:lXexplicit}, we get: $$2l_X=\pi^*\beta_*2l_Y=\pi^*\pi_*L=2L$$ as wanted.\endproof
\end{cor}
The last preliminary computation is a description of the Picard group of the generic $Y$. \begin{lemma}\label{lemma:generalPicMT}
    The generic orbifold $Y$ from the construction of \cref{ex:MTConstruction} has Picard rank $\Picrk(Y)= 2$. In particular, a sublattice of index 2 is generated by the classes $H_Y-\delta_Y,\Sigma_Y$, where $H_Y$ is the polarisation of $Y$ induced by the K3 surface $S$.
    \proof
    Recall that $S$ is defined as the double double cover $S\xrightarrow[]{\rho}P\xrightarrow[]{\mu}\PP^2$ where $\mu$ is ramified over a quartic and $\rho$ over a section of $-2K_P\cong \mu^*\mathcal{O}_{\PP^2}(2)$. A quartic in $\PP^2$ deforms in $6$ dimensions. Once fixed the quartic, since a general quartic has finite automorphism group, one can deform over the 
    family of sections of $-2K_P$.

    To compute the dimension of this family, we can describe $P$ inside the weighted projective space $\PP(1,1,1,2)$ as a section of degree 4. Notice that in this way, the canonical bundle satisfies $K_P=\mathcal{O}_{\PP_{1112}}(-1)$. Therefore, by means of the standard exact sequence
    $$0\rightarrow\mathcal{O}(4)\rightarrow \mathcal{O}_{\PP_{1112}}\rightarrow\mathcal{O}_P\rightarrow0$$
    one can conclude that $$\dim H^ 0(P,\mathcal{O}_P(2))=\dim H^ 0(\PP_{1112},\mathcal{O}_{\PP_{1112}}(2))=7$$
    as   $H^1(\PP_{1112},\mathcal{O}(8))=0$ by \cite[Section 1.4]{DolgachevWeighted}.
    From this, it follows that the Picard rank of a general $S$ is at most 8. On the other hand, as $S^{[2]}$ has a symplectic involution and it is projective $9\leq\Picrk (S^{[2]})=\Picrk (S) + 1$ implies $\Picrk (S)=8$ for general $S$, which in turn implies \[
        \Picrk\left(\faktor{X}{\sigma}\right)=1 \quad \text{and} \quad \Picrk Y=2.\]

    From \cite[Proposition 4.1]{OGradyInvolutions}, the class $H_X-\delta_X$ is invariant under the Beauville involution $\iota_0$. Since $H=\alpha^*\mathcal{O}(1)$, it is invariant for the cover involution $\tau$. Therefore $H_X$ is also invariant for $\tau^{[2]}$ and so $H_X-\delta_X$ is invariant for both $\tau^{[2]}$ and $\iota_0$, meaning it is invariant for $\sigma$.

    Therefore, $\Z(H_Y-\delta_Y)\oplus \Z \Sigma_Y$ is a finite index sublattice of $\op{Pic}Y$. In particular, since $H_X-\delta_X$ is primitive and $\eta(H^2(S^{[2]},\Z))\oplus \Z \Sigma_Y$ is an index 2 sublattice inside of $H^2(Y,\Z)$, we deduce that \[
        \Z(H_Y-\delta_Y)\oplus \Z \Sigma_Y\subset \op{Pic}Y
    \]has index 2.
    \endproof
\end{lemma}
\begin{thm}\label{lemma:classMT}
    The lagrangian fibration from the Markushevich-Tikhomirov construction (see \ref{ex:MTConstruction}) is of type A.
    \proof

    By the argument in \cite[Section 27.1]{CYmanifoldsAndRelatedGeometry}, the pullback along the flop $(\Phi\circ \beta)^*$ is an isometry. It is therefore the same to classify the monodromy orbit of $l_\cal P$ or $l_Y=(\Phi\circ \beta)^*l_\cal P$.  By \cref{lemma:discriminateThirdOrbit}, it is enough to compute the product $(l_Y,\Sigma_Y) \mod 4$.

    From \cref{cor:lYExplicit}, one can compute $(l_Y,\Sigma_Y)=\frac{1}{2}(2l_Y,\Sigma_Y)=-2k$.

    Moreover, $(l_Y,l_Y)=(l_{\cal P},l_{\cal P})=0$, so $0=(2l_Y,2l_Y)=(\beta^*\pi_*l_X,\beta^*\pi_*l_X)-4k^2=2(l_X,l_X)-4k^2$ which gives $$k^2=\frac{1}{2}(l_X,l_X).$$

    We can assume by \cref{lemma:generalPicMT} that there exist $a,k\in \Z$ such that $2l_Y=a(H_Y-\delta_Y)+k\Sigma$. As the map $\beta^*\pi_*:H^2(S^{[2]},\Z)^\sigma\rightarrow H^2(Y,\Z)$ is injective, we deduce $l_X=a(H_X-\delta_X)$ and $$k^2=\frac{1}{2}(l_X,l_X)=\frac{1}{2}a^2(H_X-\delta_X,H_X-\delta_X)=a^2$$.

    We compute now the coefficient $a$. By the Fujiki relations, $$l_X.H_X.H_X.H_X=3(l_X,H_X)(H_X,H_X)=3(aH_X,H_X)(H_X,H_X)=3\cdot 4\cdot 4a.$$ To make the notation more compact, let us denote \[
        \eta_X=\tilde \beta\circ\tilde\theta_Y:N\longrightarrow X \quad \text{and}  \quad \eta_\PP=\tilde\theta_{\cal P}\circ \phi:N\longrightarrow\dual{\PP^2}.\]

    Consider the divisor $D:=\tilde\theta_{\cal P}l _{\cal P}=\tilde\theta_{\cal P}\phi^*\mathcal{O}(1)$.
    By definition, $l_X={\eta_X}_*D$. We can now use projection formula to compute the left side of the Fujiki relations:
    $$ \int l_X.H_X.H_X.H_X=\int D.(\eta_X^*H_X^3).$$
    In the following we give a geometric description of this intersection.

    As a hyperplane section of $\dual{\PP^2}$ is described as the pencil through a chosen point $\mathbf{a}\in \PP^2$, a section of $D$ can be described as $$D_\mathbf{a}=\{(\xi,\ell)\in N, \mathbf a\in \ell\}.$$
    A section for the divisor $H_X$ is described as $$H_X:=\{\xi\in S^{[2]}, \op{supp}\xi\cap C\neq \emptyset\}$$
    where $H$ is the polarisation of $S$ and $C$ is a hyperplane section of $S$. In particular, one can deduce that a section of $(\eta_X^*H_X).^3$ is made of length 2 subschemes with support intersecting three generic planes. But 3 general hyperplane sections of $S$ will not have common points. Denoting $C_1,C_2,C_3$ the hyperplane sections, the only option for a point $\xi\in S^{[2]}$ to be supported on all three of them is to be a couple of distinct simple points $\xi=\{p,q\}$ with $p\in C_1\cap C_2, q\in C_3$. Moreover, as $S$ is degree 4 polarised, $C_1\cap C_2= 4pts$. For each one of those 4 points, moving the point $q\in C_3$ we get a curve ${\cal C}_p\subset S^{[2]}$. We get therefore 4 curves, and permuting the roles of $C_1,C_2,C_3$ we get a total of 12 curves of this kind. The intersection $H_X.^3$ is therefore the union of these 12 curves.

    Finally, we need to intersect $D_\mathbf{a}\cap{\mathcal C}_p$. Any point $(\xi,\ell)\in{\mathcal C}_p$ has to verify that $p\in \op{supp}\xi$, so $\alpha(p)\in \ell$. But then, for a point in the intersection $(\xi,\ell)\in D_\mathbf{a}\cap {\mathcal C}_p$, $\ell$ must be the unique line through $\mathbf{a}$ and $\alpha(p)$.
    In particular, if $\xi=\{p,q\}$, also $\alpha(q)$ must lie on $\ell$. Meaning that $\alpha(q)\in \alpha(C_3)\cap \ell$.
    Recall that the projection alpha can be realised as a projection from a point $\tau_\infty\in\PP^3$. As $C_3$ is a generic hyperplane section, the corresponding hyperplane in the $\PP^3$ containing $S$, doesn't pass through $\tau_\infty$. This means that $\alpha$ is injective and preserves the degree, meaning $\alpha(C_3)$ is a curve of degree 4 in $\PP^2$. Therefore there are exactly 4 choices for $q$, meaning that the intersection \[
        D_\mathbf{a}\cap{\mathcal C}_p=\{(\{p,q\},\ell), q\in \alpha^{-1}(\ell)\cap C_p\}\]
    consists of exactly four points. We therefore deduce $$\int l_X.H_X.H_X.H_X=12 |D_\mathbf{a}\cap {\cal C}_p|=12\cdot4.$$
    This, together with the Fujiki relation from before, gives
    \[
        3\cdot4\cdot 4a=\int l_X.H_X.H_X.H_X= 4\cdot 12\]
    hence $a=1$. Therefore, $k=\pm1$ and $(l_Y,\Sigma_Y)=\pm2$. By \cref{lemma:discriminateThirdOrbit} we conclude.\endproof
\end{thm}
\subsection{Consequences of the classification}
The main classification result now follows from the previous sections.
\begin{thm}\label{thm:classificationFinal}
    There are 2 deformation families of lagrangian fibrations on Nikulin-type orbifolds. In particular, let $Y$ be a Nikulin-type orbifold and $\phi: Y\rightarrow B$ a lagrangian fibration, then $\phi$ can be deformed to one of the following examples:
    \begin{itemize}
        \item A fibration induced by a lagrangian fibration defined over a $K3^{[2]}$ fourfold invariant for a symplectic involution.
        \item  A fibration induced by a lagrangian fibration defined over a $K3^{[2]}$ fourfold anti-invariant for a symplectic involution.
    \end{itemize}
    \proof
    Let $Y_i\rightarrow B_i$ for $i\in\{1,2\}$ be lagrangian fibrations induced by a lagrangian fibration defined over a $K3^{[2]}$ fourfold invariant (resp. anti-invariant) for a symplectic involution. We have shown in \cref{lemma:invariant,lemma:antiInv} that there are markings $\hat\eta_i:H^2(Y_i,\Z)\rightarrow \LambdaY$ such that $\hat\eta_i(l_i)\in\LambdaY$  are the following primitive isotropic vectors: respectively $L_1+e_2$, $L_0$. We will denote these classes respectively $\alpha_1,\alpha_2\in \LambdaY$.

    Let $\phi:Y\rightarrow B$ be a Lagrangian fibration on a Nikulin-type orbifold, $l$ the associated class, and $\eta:H^2(Y,\Z)\rightarrow \Lambda_Y$ a marking.  Write $l=k\hat l$, where $\hat l$ is primitive. By \cref{lemma:max3Orbits}, there exists a parallel transport operator $P_{\omega_2}\in Mon^2(Y)$  such that $\eta(P_{\omega_2}\hat{l})=\alpha_i$ for one $i$. Take a path $\omega_3$ in a deformation of $Y$ to $Y_i$ so that the parallel transport operator $P_{\omega_3}:H^2(Y,\Z)\rightarrow H^2(Y_i,\Z)$ induces a marking $\eta_i=\eta\circ P_{\omega_3}^{-1}:H^2(Y_i,\Z)\rightarrow \LambdaY$. By the same discussion as in \cref{rmk:anyMarking}, there exists a monodromy operator $P_{\omega_1}\in Mon^2(Y_i)$ such that $\eta_i(P_{\omega_1}l_i)=\alpha_i$. One has $\hat l= P_{\omega_2}^{-1}\eta^{-1}(\alpha_i)= P_{\omega_2}^{-1}\eta^{-1}(\eta_i(P_{\omega_1}l_i))= P_{\omega_2}^{-1}P_{\omega_3}P_{\omega_1}l_i$.  Then \cite[Theorem 3.1]{OnoratiOrtizSYZsingModuli} implies that there is a lagrangian fibration $\hat\phi:Y\rightarrow \hat B$ with associated class $m\hat l$ for some $m\in \Z_{>0}$ which is deformation equivalent to the fibration $Y_i\rightarrow B_i$. By \cref{lemma:max3Orbits} we conclude that $l=m\hat l$ and $\phi=\hat\phi$, concluding the proof.
    \endproof
\end{thm}

Notice that the polarisation type classifies the deformation type of the fibration:
\begin{cor}\label{cor:orbitsPolarisations}
    On Nikulin-type orbifold, there is a unique deformation class of lagrangian fibrations of polarisation type $(1,1)$ and a unique monodromy orbits of polarisation type $(1,2)$.
    \proof
    This simply follows from \cref{thm:classificationFinal} by \cref{rmk:polarisationTypesEquivariant}.\endproof

\end{cor}

An immediate consequence of the above classification is that the SYZ-conjecture holds for Nikulin-type orbifolds.
\begin{cor}\label{cor:SYZ}
    Let $Y$ be a Nikulin-type orbifold, let $l\in H^{1,1}(Y)$ an algebraic, isotropic, nef class. Then it is the class of a semiample line bundle. In particular, there exists a lagrangian fibration $\phi: Y\rightarrow B$ for which $ml=\phi^*\mathcal{O}_B(1)$ for some $m\in \Z_{>0}$.
    \proof
    It follows from the same proof of \cref{thm:classificationFinal} starting directly from the isotropic class $l$ instead of the lagrangian fibration $\phi$. Assume $l$ is primitive for ease of notation. With the same procedure as in the proof of the theorem, one constructs a parallel transport operator $P_\omega:H^2(Y,\Z)\rightarrow H^2(Y_i,\Z)$ for some $i$ such that $P_\omega l=l_i$ and concludes by \cite[Theorem 3.1]{OnoratiOrtizSYZsingModuli}.
    \endproof
\end{cor}

The SYZ conjecture has in turn implications for the metric geometry of the orbifold.
\begin{defin}
    The Kobayashi pseudometric on a complex variety is the maximal pseudometric such that any holomorphic map from the Poincaré disk to the variety is distance decreasing. A variety is called Kobayashi hyperbolic if the Kobayashi pseudometric is non-degenerate.
\end{defin}
In \cite[ Problem F.2, p. 405]{Kobayashi}, Kobayashi conjectured that all Calabi-Yau (meaning, canonically trivial) varieties, and so in particular IHS varieties, the Kobayashi pseudometric vanishes identically. Kamenova and Lehn proved in \cite{LehnKamenovaNonHyp} that for primitively symplectic varieties with second Betti number $b_2\geq 7$ this would follow from the SYZ conjecture. So as a consequence of \cref{cor:SYZ} and \cite[Theorem 1.1]{LehnKamenovaNonHyp} we have:
\begin{cor}
    The Kobayashi pseudometric on a Nikulin-type orbifold vanishes identically.
\end{cor}
\printbibliography
\end{document}